\documentclass{amsart}

\usepackage{url,amsmath,amssymb,amsthm,graphicx,enumitem,comment,tikz,pgfplots,longtable,multicol}

\pgfplotsset{compat=1.3}
\usepackage{tikz-network}
\usepackage[dvipsnames]{xcolor}
\usepackage[capposition=bottom]{floatrow}

\newtheorem{lemma}{Lemma}
\newtheorem{theorem}{Theorem}

\theoremstyle{definition}
\newtheorem{defn}{Definition}

\begin{document}             

\title{On the Existence of Balanced Chain Rule Task Sets}
\author[H. Gilroy]{Haile Gilroy}
\date{\today}
\address{Department of Mathematical Sciences, McNeese State University, Box 92340, Lake Charles, LA 70609, USA}
\address{Department of Mathematics \& Statistics, Auburn University, 221 Parker Hall, Auburn, AL 36849, USA} 
\email{hgilroy@mcneese.edu}
\keywords{combinatorial designs, graph decompositions, task design, Calculus I, mathematics education}
\thanks{Thank you to my colleagues, Owen Henderschedt and Jared DeLeo, for proofreading this work prior to submission.}

\let\thefootnote\relax
\footnotetext{MSC2020: Primary 05B30, Secondary 97I40.}

\begin{abstract}
In mathematics education research, mathematics task sets involving \emph{mixed practice} include tasks from many different topics within the same assignment. In this paper, we use graph decompositions to construct mixed practice task sets for Calculus I, focusing on \emph{derivative computation tasks}, or tasks of the form ``Compute $f'(x)$ of the function $f(x)=$ [elementary function]." A \emph{decomposition} $D$ of a graph $G=(V,E)$ is a collection $\{H_1, H_2, \hdots, H_t\}$ of nonempty subgraphs such that $H_i=G[E_i]$ for some nonempty subset $E_i$ of $E(G)$, and $\{E_1, E_2, \hdots, E_t\}$ is a partition of $E(G)$. We extend results on decompositions of the complete directed graph due to Meszka \& Skupie\'n to construct balanced task sets that assess the Chain Rule.
\end{abstract}

\maketitle 

\section{Introduction}

From humble rumblings about arranging schoolgirls, magic squares, and experimental design, the field of Combinatorial Design Theory has evolved into a vast area of inquiry which, over the years, has borrowed and blended ideas from many other research areas \cite{handbookchapter2}. Research in Undergraduate Mathematics Education (RUME) is an analogous field with its own well-established methodologies and research history. Graphs are particularly suited for answering structural questions about real-life scenarios \cite{ChartrandBook}, and in recent years, graphs have emerged as approaches to answering various questions in RUME. Although the use of \emph{graph-theoretic methods} in math education research has increased since 2015, this remains a minority approach among education researchers' choices of methodology. Despite this, the notion of research being an interconnected network of relationships is at the heart of many mathematics education research questions \cite{Gilroy24}. Given this promising connection between discrete mathematics and RUME and the fact that undergraduate mathematics courses, especially freshman courses such as Calculus, feature structured content curricula, it is only natural to wonder if theoretical results from Design Theory, which is both intimately related to Graph Theory and devoted to studying structural phenomena, might apply to RUME.

\begin{figure}
    \centering
    \includegraphics[width=0.5\linewidth]{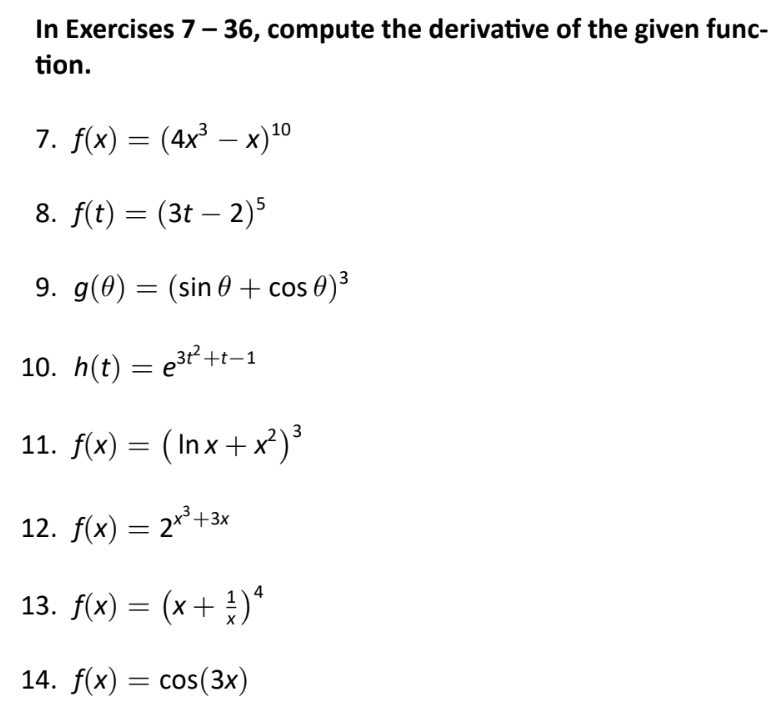}
    \caption{DC Tasks in a Calculus I textbook \cite{APEXCalc}.}
    \label{fig:DCTaskStewart}
\end{figure}

In mathematics education research, a \emph{task} is a broad term encompassing any unit of questioning assigned to a learner by an instructor. In this paper, ``task" is preferred over ``problem" or ``exercise" because education research distinguishes between the latter two depending on the level of possible mimicry of instructor-led examples by the learner \cite{Schoenfeld}. A task is \emph{skills-based} if it requires fluency in narrowly defined procedures, called \emph{skills}. A prime example of a large class of skills-based tasks is Derivative Computation tasks in Calculus I. \emph{Derivative Computation tasks} (or \emph{DC tasks}) are tasks with instructions of the form, ``Compute $f'(x)$ of the function $f(x)=$ [elementary function]," as exemplified in Figure \ref{fig:DCTaskStewart}.

Note that in this paper, we consider single-variable functions that are neither implicitly-defined functions nor functions defined by exponentiation of non-constant elementary functions.

A skills-based \emph{task set} is a collection of skills-based tasks, which, in practice, an instructor might assign to students as a homework assignment or quiz. Skills-based task sets are integral to undergraduate (especially freshman) mathematics courses. For Calculus I in particular, approximately 20-25\% of the course is devoted to computing derivatives of elementary functions using derivative rules according to a widely accepted Calculus curriculum in the U.S. \cite{AP}.

A quantitative analysis of Calculus textbooks' DC tasks has revealed that these task sets emphasize the derivatives of power functions more than any other simple function class (e.g., trigonometric functions, exponential functions, logarithmic functions, and inverse trigonometric functions) \cite{Gilroy2024b}. However, this finding misaligns with many mathematicians' perspective that Calculus students should be equally prepared for any DC skills they might encounter on assessments. Different types of misalignment in mathematics courses (including this example between instructor expectations and assigned practice), particularly in first-semester Calculus, have been linked to student retention issues in STEM majors \cite{Weston}. Thus, this paper explores the existence of what we call \emph{balanced} DC task sets, or task sets in which each DC skill in a pre-defined list occurs exactly the same number of times. In particular, the Chain Rule is notoriously difficult to teach, as described by this quote from Gordon \cite{Gordon2005}.
\begin{quotation}
    The chain rule is one of the hardest ideas to convey to students in Calculus I. It is difficult to motivate, so that most students do not really see where it comes from; it is difficult to express in symbols even after it is developed; and it is awkward to put it into words, so that many students can not remember it and so can not apply it correctly.
\end{quotation}
So, we focus our attention on constructing task sets that give students practice with this skill.

Within the research area of task design, mathematics education researchers differentiate between types of practice. \emph{Blocked practice} refers to a set of tasks that assesses one particular skill, whereas \emph{mixed practice} refers to a set of tasks that assesses many different skills \cite{BlockMixedPractice}. However, these two categories of practice are not mutually exclusive and depend on the definition of the skill being assessed. For example, when Calculus students first learn how to take derivatives, they learn some facts like the derivatives of the six trigonometric functions. These facts could constitute blocked or mixed practice depending on how one defines a skill. If we want to assess the skill of ``Taking derivatives of trigonometric functions", then recalling these facts is blocked practice. However, if we consider each fact as a skill, then recalling these facts is mixed practice. In this paper, we will consider decompositions of a vertex labeled complete directed graph because we aim to construct different types of mixed practice where the skills are the set of vertex labels. Using this approach, every possible ordered pair of skills occurs exactly once in a task set.

The remainder of this paper has three main parts. First, we establish the background of the problem, defining an index of notation in Section \ref{sec:index}, preliminary definitions in Section \ref{sec:prelims}, and known results in Section \ref{sec:priorresults}. Second, in Sections \ref{sec:NecConds} and \ref{sec:SuffConds} we extend results of Meszka \& Skupie\'n \cite{Meszka} on non-Hamiltonian directed path decompositions of the complete directed graph by proving existence results under additional conditions. Third, in Section \ref{sec:HWsamples} we give example Chain Rule task sets using our constructions, and, in Section \ref{sec:conc}, we explore possibilities for further research.

\section{Index of Notation}\label{sec:index}

\def\arraystretch{1.5}
\begin{longtable}{p{0.35\linewidth}p{0.6\linewidth}}
$\mathbb{P}$ & positive integers \\
deg$^-(v)$ & indegree of a vertex $v$\\
deg$^+(v)$ & outdegree of a vertex $v$\\
$G[E_i]$ & the subgraph of $G$ induced by the edge set $E_i\subseteq E$ \\
$\mathcal{D}K_n$ & the complete directed graph on $n$ vertices\\
DCTS($n,t$) & a derivative computation task set of order $n$ and size $t$\\
BDCTS($n,t$) & a balanced derivative computation task set of order $n$ and size $t$\\
DPD($G$) & a directed path decomposition of a directed graph $G$\\ 
HDPD($n$) & a Hamiltonian directed path decomposition of $\mathcal{D}K_n$\\
NHDPD($n$) & a non-Hamiltonian directed path decomposition of $\mathcal{D}K_n$\\
BNHDPD($n,k$) & a balanced non-Hamiltonian directed path decomposition of $\mathcal{D}K_n$ with each vertex appearing $k$ times in the decomposition\\
\end{longtable}

\section{Preliminaries}\label{sec:prelims}

In this section, we introduce terminology from Combinatorial Design Theory necessary for the results in Sections \ref{sec:NecConds} and \ref{sec:SuffConds}.\medskip

\begin{defn}[Graph Decomposition] A \emph{decomposition} $D$ of a (directed) graph $G=(V,E)$ is a collection $\{H_1, H_2, \hdots, H_t\}$ of nonempty subgraphs such that $H_i=G[E_i]$ for some nonempty subset $E_i$ of $E(G)$, and $\{E_1, E_2, \hdots, E_t\}$ is a partition of $E(G)$.\end{defn} \medskip

Figure \ref{fig:DecompExample} illustrates an example of such a decomposition. We label the vertices $v_1,\hdots, v_5$, represent the edge partition using a different color for each part, and omit single-edge parts from the picture.
\begin{figure}[h]
    \centering
    \begin{tikzpicture}[scale=1]
        \Vertex[x=0,y=2,color=white,label=$v_1$]{1}
        \Vertex[x=1.902,y=0.618,color=white,label=$v_2$]{2}
        \Vertex[x=1.176,y=-1.618,color=white,label=$v_3$]{3}
        \Vertex[x=-1.176,y=-1.618,color=white,label=$v_4$]{4}
        \Vertex[x=-1.902,y=0.618,color=white,label=$v_5$]{5}

    \Edge[color=blue,Direct,bend=10](4)(5)
    \Edge[color=blue,Direct,bend=10](5)(1)
    \Edge[color=blue,Direct,bend=10](1)(2)

    \Edge[color=red,Direct,bend=10](1)(5)
    \Edge[color=red,Direct,bend=10](5)(4)
    \Edge[color=red,Direct,bend=10](4)(3)

    \Edge[color=yellow,Direct,bend=10](2)(3)
    \Edge[color=yellow,Direct,bend=10](3)(4)
    \Edge[color=yellow,Direct,bend=10](4)(1)

    \Edge[color=green,Direct,bend=10](3)(2)
    \Edge[color=green,Direct,bend=10](2)(1)

    \Edge[color=cyan,Direct,bend=10](3)(1)
    \Edge[color=cyan,Direct,bend=10](1)(4)

    \Edge[color=orange,Direct,bend=10](4)(2)
    \Edge[color=orange,Direct,bend=10](2)(5)

    \Edge[color=pink,Direct,bend=10](1)(3)
    \Edge[color=pink,Direct,bend=10](3)(5)
    \end{tikzpicture}
    \caption{A decomposition of $\mathcal{D}K_5$ into directed paths of varying lengths}   
    \label{fig:DecompExample}
\end{figure}

Many graph decomposition results concern unlabeled graphs. However, when constructing Derivative Computation task sets, it makes sense to consider decompositions of (vertex) labeled graphs so we can specify the simple functions and arithmetic operations present in an elementary function. 

\begin{defn}[Elementary Function Vertex Labeling]
    An \emph{elementary function vertex labeling} of a graph $G$ is a (not necessarily injective) function $\Lambda: V\to\mathcal{L}$ where $\mathcal{L}=\mathcal{F}\cup\mathcal{O}$ is a set of labels such that:
\begin{enumerate}
    \item $\mathcal{F}$ is a set of simple function classes, and
    \item $\mathcal{O}\subseteq\{+_i\ |\ i\in\mathbb{P}, i\geq2\}\cup\{\times_i\ |\ i\in\mathbb{P}, i\geq2\}\cup\{\div\}$ is a set of operations where $i$ denotes the number of summands in a sum or factors in a product, respectively.
\end{enumerate}
\end{defn}

In the results that follow, we decompose $\mathcal{D}K_n$ into a particular class of directed trees, which we call elementary function trees. Similar representations of functions, called \emph{syntax trees}, are a common data structure in Computer Science. 

\begin{defn}[Construction Version of EFT]\label{defn:EFTConst}
    An \emph{elementary function tree}, or EFT, is a labeled, rooted in-tree, $(V,E,\Lambda)$, where 
    \begin{enumerate}
        \item $V$ is a collection of simple functions and addition, multiplication, and division operations
        \item specified by a labeling $\Lambda$, and
        \item $E$ is the set of composition relations among the vertex labels (i.e. for $u,v\in V, uv\in E$ iff $\Lambda(v)\circ\Lambda(u)$).
    \end{enumerate}
\end{defn}

A natural characterization of EFTs exists based on whether they correspond to practical tasks or not.
\begin{enumerate}
    \item A \emph{feasible} EFT is an EFT that has at least one corresponding elementary function. In other words, if $v\in V$ with $\Lambda(v)\in\mathcal{F}$, then deg$^-(v)\leq1$ and deg$^+(v)\leq1$, and if $v\in V$ with $\Lambda(v)\in\mathcal{O}$, then deg$^-(v)=i$ and deg$^+(v)\leq1$.
    \item A \emph{semi-feasible} EFT is an EFT that can be made feasible by augmenting it with a finite number of additional vertices and arcs. In other words, $\exists\ v\in V$ with $\Lambda(v)\in\mathcal{O}$ such that deg$^-(v)<i$.
    \item An \emph{infeasible} EFT is an EFT that does not have a corresponding elementary function. In other words, either $\exists\ v\in V$ with $\Lambda(v)\in\mathcal{F}$ such that either deg$^-(v)>1$ or deg$^+(v)>1$, or $\exists\ v\in V$ with $\Lambda(v)\in\mathcal{O}$ such that deg$^+(v)>1$.
\end{enumerate}
Figure \ref{fig:EFTClasses} depicts an example of each class of EFTs. The leftmost EFT is feasible because we can find an elementary function it represents; for example, consider the function $f(x)=(x^2+\sin x)^2$. The middle EFT is semi-feasible because we can make it feasible by adding two simple function vertices pointing into the multiplication vertex. The rightmost EFT is infeasible because a single-variable trigonometric function cannot have two input values. However, note that if we were to consider multivariable functions, these notions of feasible and infeasible would change.

\begin{figure}[h]
\centering
\begin{multicols}{3}
\begin{tikzpicture}[scale=0.75]
    \draw[] (-2,0) circle (.5cm) node{$P$};
    \draw[<-] (-1.5,0)--(-0.5,0);
    \draw[] (0,0) circle (.5cm) node{$+_2$};
    \draw[<-,rotate=-45] (.5,0)--(1.5,0);
    \draw[<-,rotate=45] (.5,0)--(1.5,0);
    \draw[] (1.45,1.44) circle (.5cm) node{$P$};
    \draw[] (1.45,-1.44) circle (.5cm) node{$T$};
    \draw[] (0,-2.5) node{(a) feasible EFT};
\end{tikzpicture}

\begin{tikzpicture}[scale=0.75]
    \draw[] (-2,0) circle (.5cm) node{$P$};
    \draw[<-] (-1.5,0)--(-0.5,0);
    \draw[] (0,0) circle (.5cm) node{$+_2$};
    \draw[<-,rotate=-45] (.5,0)--(1.5,0);
    \draw[<-,rotate=45] (.5,0)--(1.5,0);
    \draw[] (1.45,1.44) circle (.5cm) node{$\times_2$};
    \draw[] (1.45,-1.44) circle (.5cm) node{$T$};
    \draw[] (0,-2.5) node{(b) semi-feasible EFT};
\end{tikzpicture}

\begin{tikzpicture}[scale=0.75]
    \draw[] (-2,0) circle (.5cm) node{$P$};
    \draw[<-] (-1.5,0)--(-0.5,0);
    \draw[] (0,0) circle (.5cm) node{$T$};
    \draw[<-,rotate=-45] (.5,0)--(1.5,0);
    \draw[<-,rotate=45] (.5,0)--(1.5,0);
    \draw[] (1.45,1.44) circle (.5cm) node{$P$};
    \draw[] (1.45,-1.44) circle (.5cm) node{$+_2$};
    \draw[] (0,-2.5) node{(c) infeasible EFT};
\end{tikzpicture}
\end{multicols}
\caption{A feasible, a semi-feasible, and an infeasible EFT. The labels $P$, $T$, $+_2$, and $\times_2$ stand for ``power function", ``trigonometric function", ``sum of two addends", and ``product of two factors", respectively.}
    \label{fig:EFTClasses}
\end{figure}

A feasible EFT is synonymous with a DC task. Therefore, we define a \emph{DC task set}, $D$, as a collection of feasible EFTs. The \emph{size} of the task set, $t=|D|$, is the number of tasks (or feasible EFTs) that it contains, whereas the \emph{order} of the task set is $n=|\mathcal{L}|=V(G)$. We denote a DC task set of order $n$ as DCTS($n$).  
If every label appears exactly the same number of times in a task set, then we say the DCTS is \emph{balanced}. We denote a balanced DC task set of order $n$ by BDCTS($n$).

\section{Known Results}\label{sec:priorresults}

This section discusses known results on directed path decompositions, DPDs, and how they may be used to construct Derivative Computation Task Sets of order $n$, or DCTS($n$). From a graph-theoretic point of view, the simplest EFTs are directed paths, which correspond to DC tasks that assess the Chain Rule in isolation. Lemma \ref{lem:DPDisDCTS} allows us to restrict our focus to constructing DPDs with special properties. 

\begin{lemma}\label{lem:DPDisDCTS}
    If $D$ is a DPD$(G)$ with $\mathcal{L}=\mathcal{F}$ and $|\mathcal{L}|=|V(G)|=n$, then $D$ is a DCTS$(n)$. 
\end{lemma}

\begin{proof}
    Suppose $D$ is a DPD($G$). Then, all EFTs in $D$ are directed paths. This means that, for every vertex $v$ contained in any directed path $P$ in $D$, deg$^-$($v$)$\leq 1$ and deg$^+$($v$)$\leq 1$. Since $\mathcal{L}=\mathcal{F}$ ($\mathcal{O}=\varnothing$) and all simple functions have exactly one input, all directed paths in $D$ are feasible. Thus, $D$ is a DCTS($n$).
\end{proof}

\subsection{Directed Path Decompositions of $\mathcal{D}K_n$}

A natural starting point for studying DPD($\mathcal{D}K_n$) is to require only Hamiltonian paths in the decomposition. The existence problem of HDPD($n$) was completely solved by Bos\'ak \cite{Bosak1990} by synthesizing and extending results of Bermond and Faber \cite{BermondFaber} and Tillson \cite{Tillson}. However, if we are to assign DCTS to students, Hamiltonian paths are impractical. Each task in a HDPD($n$) would require students to practice the Chain Rule $n-1$ times, meaning we would need to construct task sets with relatively small values of $n$, restricting the variety in simple function classes students would be exposed to in a task set.

Considering the impracticality of Hamiltonian paths, the next natural question is to consider non-Hamiltonian paths, or NHDPD($n$). A theorem of Meszka and Skupie\'n \cite{Meszka} completely solves the existence problem of NHDPD($n$). 

\begin{theorem}[M. Meszka \& Z. Skupie\'n, 2006]\label{thm:MSNHPaths}
    NHDPD($n$) with paths of arbitrarily prescribed lengths $(\leq n-2)$ exist for any positive integer $n\geq3$, provided that the lengths sum up to $n(n-1)$, the size of $\mathcal{D}K_n$.
\end{theorem}

We extend the notion of NHDPD($n$) by considering their balanced counterparts.

\begin{defn}
    Let $n,k\in\mathbb{P}$. A \emph{balanced non-Hamiltonian directed path decomposition} of $\mathcal{D}K_n$, denoted BNHDPD($n,k$), is a NHDPD($n$) such that each vertex is contained in exactly $k$ paths in the decomposition.
\end{defn}

By Lemma \ref{lem:DPDisDCTS}, all BNHDPD($n,k$) are BDCTS($n$) where each elementary function occurs $k$ times in the task set. \medskip

\section{Necessary Conditions for BNHDPD($n,k$)}\label{sec:NecConds}

We now give general necessary conditions for the existence of a BNHDPD$(n,k)$.

\begin{lemma}
    If a BNHDPD$(n,k)$ comprised of $x_i$ non-Hamiltonian directed paths of length $i$ exists, then \begin{equation}\label{eq:nhpathsnec}
        \sum_{i=1}^{n-2}ix_i=n(n-1).
    \end{equation}
\end{lemma}

\begin{proof}
    Since such a decomposition is a non-Hamiltonian directed path decomposition of the complete directed graph, this follows directly from Theorem \ref{thm:MSNHPaths}.
\end{proof}

The \emph{incidence matrix} of a BNHDPD$(n,k)$, $A$, is an array whose rows represent the vertices of $\mathcal{D}K_n$, $v_1,\hdots,v_n$, and whose columns represent the EFTs of the task set, $T_1, \hdots, T_{|D|}$ such that
\[[A]=\begin{cases}
    1\ \text{if}\ v_i\in T_j \\
    0\ \text{if}\ v_i\notin T_j
\end{cases}\]

\begin{lemma}
    If a BNHDPD$(n,k)$ comprised of $x_i$ non-Hamiltonian directed paths of length $i$ exists, then \begin{equation}\label{eq:balnhpathsnec}
    \sum_{i=1}^{n-2}(i+1)x_i=nk,
    \end{equation}
    where $k$ is the number of paths incident to each vertex.
\end{lemma}

\begin{proof}
    We count the total number of ones in the incidence matrix of a BNHDPD$(n,k)$ in two different ways. First, count by rows. Since there are $n$ vertices, each incident to $k$ paths, there are $nk$ ones in the incidence matrix. Now, count by columns. Since there are $x_i$ paths of length $i$ and $(i+1)$ vertices in a path of length $i$, there are $2x_1+3x_2+\cdots+(n-1)x_{n-2}$ ones in the incidence matrix. This gives the desired result.
\end{proof}

\begin{theorem}\label{Thm:BalDCTS:Nec} If a BNHDPD$(n,k)$ comprised of $x_i$ non-Hamiltonian directed paths of length $i$ exists, then $n$ divides the number of paths in the decomposition. \end{theorem}

\begin{proof}
    Substituting (\ref{eq:nhpathsnec}) into (\ref{eq:balnhpathsnec}), we arrive at 
    \begin{equation}\sum_{i=1}^{n-2}x_i=n(k-(n-1)),\end{equation}\label{eq:balnhpathsnec2}
    which is the desired result. 
\end{proof}

The next necessary condition concerns the structure of a BNHDPD$(n,k)$ at each vertex of $\mathcal{D}K_n$. 

\begin{theorem}\label{Thm:BalDCTS:Nec2}
If a BNHDPD$(n,k)$ comprised of $x_i$ non-Hamiltonian directed paths of length $i$ exists $(i\in\mathbb{P},\leq n-2)$, then $n$ divides the number of interior vertices among all paths in the decomposition. That is,
\[n\ \biggr|\ \sum_{i=2}^{n-2} (i-1)x_i\]
\end{theorem}

\begin{proof}
    \begin{align} 
      &\sum_{i=2}^{n-2} (i-1)x_i \nonumber\\
    = &\sum_{i=2}^{n-2} ix_i-\sum_{i=2}^{n-2} x_i
    \end{align}\label{eq:balnhpathsnec3}
    Substituting (\ref{eq:nhpathsnec}) and (3) into (4), we arrive at 
    \[\sum_{i=2}^{n-2} (i-1)x_i=n[n-1-(k-(n-1))],\] which is the desired result.
    
\end{proof}

\section{Sufficient Conditions}\label{sec:SuffConds}

In Subsections \ref{ssec:BNHDPD:n=5} and \ref{ssec:BNHDPD:n=6}, we establish sufficiency for $n=5$ and $n=6$, respectively. We chose to focus on these cases because they produce a variety of practical, balanced task sets to assign to students. The case $n=4$ results in little possible variety (see Figure \ref{fig:Task Set Enum}), and only two of these possible task sets meet the divisibility condition for balance. The case $n=7$ is impractical for two reasons. First, paths of length 5 are introduced, meaning each task can require up to 5 successive applications of the chain rule, which we believe is unnecessarily complex for students. Second, though $n=7$ has greater potential variety in tasks, the task sets will be larger than what we would assign to students.

Theorem \ref{thm:MSNHPaths} allows us to represent any task set as an nonnegative integer solution to the Diophantine Equation 
\[\sum_{i=1}^{n-2} ix_i=n(n-1).\]
Figure \ref{fig:Task Set Enum} shows the number of nonnegative integer solutions to this equation for $n=4,5,6$, categorized by their sum (the size of the task set). This gives us the maximum number of sufficiency cases we will need to construct in Subsections \ref{ssec:BNHDPD:n=5} and \ref{ssec:BNHDPD:n=6} and hints at the wide variety in possible task sets using our constructions.

\begin{figure}[h]
\centering  

\begin{tikzpicture}[scale=1,trim axis left,trim axis right]
\begin{axis}[
    ybar,
    height=3in,
    width=5in,
    xtick distance=1,
	xticklabel style={
		text width=7cm, align=center, font=\footnotesize},
	ylabel=Frequency,
    xlabel=$\sum x_i$,
    ymin=0,
    xmin=5,
    xmax=31,
    bar width=10pt,
    title={Enumeration of Some Nonnegative Integer Solutions to $\sum ix_i=n(n-1)$},
	legend style={at={(0.85,0.95)},
	anchor=north,legend columns=1},
]
\addplot[black, fill=gray, fill opacity=0.5, bar shift=-0.5]
    coordinates {(8,2) (9,7) (10,14) (11,20) (12,25) (13,26) (14,28) (15,27) (16,24) (17,20) (18,18) (19,16) (20,14) (21,12) (22,10) (23,8) (24,7) (25,5) (26,4) (27,3) (28,2) (29,1) (30,1)}; \addlegendentry[text width=45pt, text depth=]{$n=6$};

\addplot[blue, fill=blue, fill opacity=0.5, bar shift=-0.5] 
	coordinates {(8,2) (9,4) (10,6) (11,5) (12,5) (13,4) (14,4) (15,3) (16,3) (17,2) (18,2) (19,1) (20,1)}; \addlegendentry[text width=45pt, text depth=]{$n=5$};

\addplot[orange, fill=orange, fill opacity=0.5,bar shift=-0.5]
    coordinates {(6,1) (7,1) (8,1) (9,1) (10,1) (11,1) (12,1)}; \addlegendentry[text width=45pt, text depth=]{$n=4$};
    
\end{axis}
\end{tikzpicture}
\caption{Enumeration of some possible DCTS$(n,t)$ to be constructed from NHDPD($\mathcal{D}K_n$).}
\label{fig:Task Set Enum}
\end{figure}

Before we establish sufficiency, we define the reverse of a directed path, which appears in many of our constructions.    Given a directed path $P=v_1v_2\cdots v_m$, the \emph{reverse} of $P$ is the directed path $P'=v_mv_{m-1}\cdots v_1$. We also establish one general sufficiency case.

\begin{lemma}\label{Lem:BNHDPD:Trivial}
    A BNHDPD$(n,k)$ with $n(n-1)$ directed paths of length 1 exists for all $n\in\mathbb{P}$, $n\geq2$.
\end{lemma}

\begin{proof}
    Since every arc of $\mathcal{D}K_n$ is a directed path of length 1, it follows that the edge set of $\mathcal{D}K_n$ is a BNHDPD$(n,2(n-1))$.
\end{proof}

\subsection{Sufficiency for BNHDPD$(5,k)$}\label{ssec:BNHDPD:n=5}

In this subsection, we establish sufficient conditions for BNHDPD$(n,k)$ with $n=5$. This is equivalent to constructing BDCTS$(n)$ that contain 5 classes of simple functions and that assess the Chain Rule in isolation. Sample task sets constructed using results in this subsection are given in Section \ref{sec:HWsamples}.

\begin{theorem}\label{Thm:BNHDPD:Suff5}
    A BNHDPD$(5,k)$ comprised of $x_1$ copies of $\mathcal{D}P_2$, $x_2$ copies of $\mathcal{D}P_3$, and $x_3$ copies of $\mathcal{D}P_4$ exists for all nonnegative integer solutions of $x_1+2x_2+3x_3=20$ such that $5\ |\ x_1+x_2+x_3$.
\end{theorem}

\begin{proof}
    Let $V(\mathcal{D}K_5)=\{v_1,\hdots, v_5\}$. We proceed by construction via case analysis. Since $5\ |\ x_1+x_2+x_3$ and $x_1+x_2+x_3\leq 20$, we only have to verify solutions of $x_1+x_2+x_3=s$ for $s\in\{10,15\}$. The case $s=20$ is already established by Lemma \ref{Lem:BNHDPD:Trivial}. Within each case, we denote each subcase by the ordered triple $(x_1,x_2,x_3)$ of the corresponding solution to $x_1+x_2+x_3=s$, for $s\in\{10,15\}$. For each subcase, we construct the BNHDPD$(5,k)$
    \[\bigcup_{i=1}^{n-2} (P_i\cup P'_i).\] where $P_i, 1\leq i\leq n-2$, is the set of directed paths of length $i$, and  $P'_i, 1\leq i\leq n-2$, is the set of some, but not necessarily all, reverses of $P_i$. Visualizations of each construction are given in Figures \ref{fig:n=5:case1} and \ref{fig:n=5:case2}.
    \medskip

    \noindent\textbf{Case 1:} $x_1+x_2+x_3=10$
   
    \begin{enumerate}[label=\textit{Subcase 1\alph*.}]
       {\setlength\itemindent{38pt} \item $(0,10,0)$}
       
       $P_2=\{v_1v_2v_4, v_2v_3v_5, v_3v_4v_1, v_4v_5v_2, v_5v_1v_3\}$

       $P'_2=\{v_4v_2v_1, v_5v_3v_2, v_1v_4v_3, v_2v_5v_4, v_3v_1v_5\}$
       
       $P_3, P'_3,P_1,P'_1=\varnothing$

       {\setlength\itemindent{38pt} \item $(1,8,1)$}
       
       $P_3=\{v_5v_1v_2v_4\}$
       
       $P_2=\{v_2v_3v_5, v_3v_4v_1, v_4v_5v_2, v_2v_1v_3, v_1v_5v_3, v_5v_4v_2, v_4v_3v_1, v_3v_2v_5\}$

       $P_1=\{v_1v_4\}$

       $P'_3,P'_2,P'_1=\varnothing$

       {\setlength\itemindent{38pt} \item $(2,6,2)$}
       
       $P_3=\{v_4v_5v_1v_2\}$

       $P'_3=\{v_2v_1v_5v_4\}$

       $P_2=\{v_2v_3v_1, v_3v_4v_1, v_4v_2v_5\}$

       $P'_2=\{v_1v_3v_2, v_1v_4v_3, v_5v_2v_4\}$

       $P_1=\{v_3v_5\}$

       $P'_1=\{v_5v_3\}$

       {\setlength\itemindent{38pt} \item $(3,4,3)$}
       
       $P_3=\{v_4v_5v_1v_2, v_1v_5v_4v_3, v_2v_3v_4v_1\}$
       
       $P_2=\{v_3v_2v_1, v_3v_1v_4, v_4v_2v_5, v_1v_3v_5\}$

       $P_1=\{v_5v_3, v_5v_2, v_2v_4\}$

       $P'_3,P'_2,P'_1=\varnothing$

       {\setlength\itemindent{38pt} \item $(4,2,4)$}
       
       $P_3=\{v_1v_2v_4v_3, v_2v_3v_5v_4, v_3v_4v_1v_5, v_4v_5v_2v_1\}$
       
       $P_2=\{v_3v_1v_4, v_5v_3v_2\}$

       $P_1=\{v_5v_1, v_1v_3, v_2v_5, v_4v_2\}$

       $P'_3,P'_2,P'_1=\varnothing$

       {\setlength\itemindent{36pt} \item $(5,0,5)$}
       
       $P_3=\{v_1v_2v_4v_3, v_2v_3v_5v_4, v_3v_4v_1v_5, v_4v_5v_2v_1, v_5v_1v_3v_2\}$

       $P_1=\{v_2v_5, v_3v_1, v_4v_2, v_5v_3, v_1v_4\}$

       $P'_3,P_2,P'_2,P'_1=\varnothing$ 
              
    \end{enumerate}

    \begin{figure}[h]
        \centering
        \begin{tabular}{|c|c|c|}
\hline
& & \\

\begin{tikzpicture}
\Vertex[x=0,y=2,color=white,label=$v_1$]{1}
\Vertex[x=1.902,y=0.618,color=white,label=$v_2$]{2}
\Vertex[x=1.176,y=-1.618,color=white,label=$v_3$]{3}
\Vertex[x=-1.176,y=-1.618,color=white,label=$v_4$]{4}
\Vertex[x=-1.902,y=0.618,color=white,label=$v_5$]{5}

\Edge[color=blue](1)(2)
\Edge[color=blue](2)(4)

\Edge[color=red](5)(1)
\Edge[color=red](1)(3)

\Edge[color=cyan](2)(3)
\Edge[color=cyan](3)(5)

\Edge[color=green](3)(4)
\Edge[color=green](4)(1)

\Edge[color=yellow](4)(5)
\Edge[color=yellow](5)(2)

\end{tikzpicture} &

\begin{tikzpicture}
\Vertex[x=0,y=2,color=white,label=$v_1$]{1}
\Vertex[x=1.902,y=0.618,color=white,label=$v_2$]{2}
\Vertex[x=1.176,y=-1.618,color=white,label=$v_3$]{3}
\Vertex[x=-1.176,y=-1.618,color=white,label=$v_4$]{4}
\Vertex[x=-1.902,y=0.618,color=white,label=$v_5$]{5}

\Edge[color=blue,Direct,bend=10](5)(1)
\Edge[color=blue,Direct,bend=10](1)(2)
\Edge[color=blue,Direct,bend=10](2)(4)

\Edge[color=red,Direct,bend=10](2)(3)
\Edge[color=red,Direct,bend=10](3)(5)

\Edge[color=yellow,Direct,bend=10](3)(4)
\Edge[color=yellow,Direct,bend=10](4)(1)

\Edge[color=green,Direct,bend=10](4)(5)
\Edge[color=green,Direct,bend=10](5)(2)

\Edge[color=orange,Direct,bend=10](2)(1)
\Edge[color=orange,Direct,bend=10](1)(3)

\Edge[color=violet,Direct,bend=10](1)(5)
\Edge[color=violet,Direct,bend=10](5)(3)

\Edge[color=cyan,Direct,bend=10](5)(4)
\Edge[color=cyan,Direct,bend=10](4)(2)

\Edge[color=pink,Direct,bend=10](4)(3)
\Edge[color=pink,Direct,bend=10](3)(1)

\Edge[color=JungleGreen,Direct,bend=10](3)(2)
\Edge[color=JungleGreen,Direct,bend=10](2)(5)

\end{tikzpicture} &

\begin{tikzpicture}
\Vertex[x=0,y=2,color=white,label=$v_1$]{1}
\Vertex[x=1.902,y=0.618,color=white,label=$v_2$]{2}
\Vertex[x=1.176,y=-1.618,color=white,label=$v_3$]{3}
\Vertex[x=-1.176,y=-1.618,color=white,label=$v_4$]{4}
\Vertex[x=-1.902,y=0.618,color=white,label=$v_5$]{5}

\Edge[color=blue](4)(5)
\Edge[color=blue](5)(1)
\Edge[color=blue](1)(2)

\Edge[color=red](2)(3)
\Edge[color=red](3)(1)

\Edge[color=yellow](3)(4)
\Edge[color=yellow](4)(1)

\Edge[color=green](4)(2)
\Edge[color=green](2)(5)

\end{tikzpicture} \\

Subcase (a) & Subcase (b) & Subcase (c) \\
\hline

& &  \\

\begin{tikzpicture}
\Vertex[x=0,y=2,color=white,label=$v_1$]{1}
\Vertex[x=1.902,y=0.618,color=white,label=$v_2$]{2}
\Vertex[x=1.176,y=-1.618,color=white,label=$v_3$]{3}
\Vertex[x=-1.176,y=-1.618,color=white,label=$v_4$]{4}
\Vertex[x=-1.902,y=0.618,color=white,label=$v_5$]{5}

\Edge[color=blue,Direct,bend=10](4)(5)
\Edge[color=blue,Direct,bend=10](5)(1)
\Edge[color=blue,Direct,bend=10](1)(2)

\Edge[color=red,Direct,bend=10](1)(5)
\Edge[color=red,Direct,bend=10](5)(4)
\Edge[color=red,Direct,bend=10](4)(3)

\Edge[color=yellow,Direct,bend=10](2)(3)
\Edge[color=yellow,Direct,bend=10](3)(4)
\Edge[color=yellow,Direct,bend=10](4)(1)

\Edge[color=green,Direct,bend=10](3)(2)
\Edge[color=green,Direct,bend=10](2)(1)

\Edge[color=cyan,Direct,bend=10](3)(1)
\Edge[color=cyan,Direct,bend=10](1)(4)

\Edge[color=orange,Direct,bend=10](4)(2)
\Edge[color=orange,Direct,bend=10](2)(5)

\Edge[color=pink,Direct,bend=10](1)(3)
\Edge[color=pink,Direct,bend=10](3)(5)

\end{tikzpicture} &

\begin{tikzpicture}
\Vertex[x=0,y=2,color=white,label=$v_1$]{1}
\Vertex[x=1.902,y=0.618,color=white,label=$v_2$]{2}
\Vertex[x=1.176,y=-1.618,color=white,label=$v_3$]{3}
\Vertex[x=-1.176,y=-1.618,color=white,label=$v_4$]{4}
\Vertex[x=-1.902,y=0.618,color=white,label=$v_5$]{5}

\Edge[color=blue,Direct,bend=10](1)(2)
\Edge[color=blue,Direct,bend=10](2)(4)
\Edge[color=blue,Direct,bend=10](4)(3)

\Edge[color=red,Direct,bend=10](2)(3)
\Edge[color=red,Direct,bend=10](3)(5)
\Edge[color=red,Direct,bend=10](5)(4)

\Edge[color=yellow,Direct,bend=10](3)(4)
\Edge[color=yellow,Direct,bend=10](4)(1)
\Edge[color=yellow,Direct,bend=-10](1)(5)

\Edge[color=green,Direct,bend=10](4)(5)
\Edge[color=green,Direct,bend=10](5)(2)
\Edge[color=green,Direct,bend=10](2)(1)

\Edge[color=cyan,Direct,bend=10](3)(1)
\Edge[color=cyan,Direct,bend=10](1)(4)

\Edge[color=pink,Direct,bend=10](5)(3)
\Edge[color=pink,Direct,bend=10](3)(2)

\end{tikzpicture} &

\begin{tikzpicture}
\Vertex[x=0,y=2,color=white,label=$v_1$]{1}
\Vertex[x=1.902,y=0.618,color=white,label=$v_2$]{2}
\Vertex[x=1.176,y=-1.618,color=white,label=$v_3$]{3}
\Vertex[x=-1.176,y=-1.618,color=white,label=$v_4$]{4}
\Vertex[x=-1.902,y=0.618,color=white,label=$v_5$]{5}

\Edge[color=red,Direct,bend=10](5)(1)
\Edge[color=red,Direct,bend=10](1)(3)
\Edge[color=red,Direct,bend=10](3)(2)

\Edge[color=blue,Direct,bend=10](1)(2)
\Edge[color=blue,Direct,bend=10](2)(4)
\Edge[color=blue,Direct,bend=10](4)(3)

\Edge[color=cyan,Direct,bend=10](2)(3)
\Edge[color=cyan,Direct,bend=10](3)(5)
\Edge[color=cyan,Direct,bend=10](5)(4)

\Edge[color=green,Direct,bend=10](3)(4)
\Edge[color=green,Direct,bend=10](4)(1)
\Edge[color=green,Direct,bend=10](1)(5)

\Edge[color=yellow,Direct,bend=10](4)(5)
\Edge[color=yellow,Direct,bend=10](5)(2)
\Edge[color=yellow,Direct,bend=10](2)(1)

\end{tikzpicture}\\

Subcase (d) & Subcase (e) & Subcase (f)\\

\hline

\end{tabular}
        \caption{Visual Constructions from Theorem \ref{Thm:BNHDPD:Suff5}, Case 1}
        \floatfoot{Each color corresponds to a different path in its corresponding construction. Undirected paths represent a directed path and its reverse. Paths of length 1 are not included in the visualizations.}
        \label{fig:n=5:case1}
    \end{figure}\newpage

    \noindent\textbf{Case 2:} $x_1+x_2+x_3=15$

    \begin{enumerate}[label=\textit{Subcase 2\alph*.}]
   
       {\setlength\itemindent{38pt} \item $(10,5,0)$}
       
       $P_2=\{v_5v_1v_2, v_2v_3v_4, v_4v_5v_2, v_5v_4v_3, v_3v_2v_1\}$

       $P_1=\{v_1v_5, v_2v_5, v_1v_3, v_1v_4, v_2v_4, v_3v_5\}$

       $P'_1=\{v_3v_1, v_4v_1, v_4v_2, v_5v_3\}$

       $P_3,P'_3,P'_2=\varnothing$

       {\setlength\itemindent{38pt} \item $(11,3,1)$}
       
       $P_3=\{v_5v_1v_2v_4\}$
       
       $P_2=\{v_2v_3v_5, v_1v_4v_2, v_2v_5v_4\}$

       $P_1=\{v_2v_1, v_1v_5, v_4v_5, v_5v_2, v_3v_2, v_5v_3, v_4v_1, v_1v_3, v_3v_4\}$

       $P'_1=\{v_3v_1, v_4v_3\}$

       $P'_3,P'_2=\varnothing$ 
       
       \newpage

       {\setlength\itemindent{38pt} \item $(12,1,2)$}
       
       $P_3=\{v_5v_1v_2v_4, v_2v_3v_4v_5\}$
       
       $P_2=\{v_3v_5v_2\}$

       $P_1=\{v_2v_5, v_2v_1, v_1v_5, v_5v_3, v_5v_4, v_4v_2, v_4v_3, v_3v_2, v_1v_3, v_1v_4\}$

       $P'_1=\{v_3v_1, v_4v_1\}$
       
       $P'_3,P'_2=\varnothing$ 
       \end{enumerate}

\begin{figure}[h]
    \centering
    \begin{tabular}{|c|c|c|}
\hline
& & \\

\begin{tikzpicture}
\Vertex[x=0,y=2,color=white,label=$v_1$]{1}
\Vertex[x=1.902,y=0.618,color=white,label=$v_2$]{2}
\Vertex[x=1.176,y=-1.618,color=white,label=$v_3$]{3}
\Vertex[x=-1.176,y=-1.618,color=white,label=$v_4$]{4}
\Vertex[x=-1.902,y=0.618,color=white,label=$v_5$]{5}

\Edge[color=blue,Direct,bend=10](5)(1)
\Edge[color=blue,Direct,bend=10](1)(2)

\Edge[color=red,Direct,bend=10](2)(3)
\Edge[color=red,Direct,bend=10](3)(4)

\Edge[color=yellow,Direct,bend=10](4)(5)
\Edge[color=yellow,Direct,bend=10](5)(2)

\Edge[color=green,Direct,bend=10](5)(4)
\Edge[color=green,Direct,bend=10](4)(3)

\Edge[color=cyan,Direct,bend=10](3)(2)
\Edge[color=cyan,Direct,bend=10](2)(1)

\end{tikzpicture} &

\begin{tikzpicture}
\Vertex[x=0,y=2,color=white,label=$v_1$]{1}
\Vertex[x=1.902,y=0.618,color=white,label=$v_2$]{2}
\Vertex[x=1.176,y=-1.618,color=white,label=$v_3$]{3}
\Vertex[x=-1.176,y=-1.618,color=white,label=$v_4$]{4}
\Vertex[x=-1.902,y=0.618,color=white,label=$v_5$]{5}

\Edge[color=blue,Direct,bend=10](5)(1)
\Edge[color=blue,Direct,bend=10](1)(2)
\Edge[color=blue,Direct,bend=10](2)(4)

\Edge[color=red,Direct,bend=10](2)(3)
\Edge[color=red,Direct,bend=10](3)(5)

\Edge[color=yellow,Direct,bend=10](1)(4)
\Edge[color=yellow,Direct,bend=10](4)(2)

\Edge[color=green,Direct,bend=10](2)(5)
\Edge[color=green,Direct,bend=-10](5)(4)

\end{tikzpicture} &

\begin{tikzpicture}
\Vertex[x=0,y=2,color=white,label=$v_1$]{1}
\Vertex[x=1.902,y=0.618,color=white,label=$v_2$]{2}
\Vertex[x=1.176,y=-1.618,color=white,label=$v_3$]{3}
\Vertex[x=-1.176,y=-1.618,color=white,label=$v_4$]{4}
\Vertex[x=-1.902,y=0.618,color=white,label=$v_5$]{5}

\Edge[color=blue,Direct,bend=10](5)(1)
\Edge[color=blue,Direct,bend=10](1)(2)
\Edge[color=blue,Direct,bend=10](2)(4)

\Edge[color=red,Direct,bend=10](2)(3)
\Edge[color=red,Direct,bend=10](3)(4)
\Edge[color=red,Direct,bend=10](4)(5)

\Edge[color=yellow,Direct,bend=10](3)(5)
\Edge[color=yellow,Direct,bend=10](5)(2)

\end{tikzpicture} \\

Subcase (a) & Subcase (b) & Subcase (c) \\
\hline

\end{tabular}
    \caption{Visual Constructions from Theorem \ref{Thm:BNHDPD:Suff5}, Case 2}
    \floatfoot{Each color corresponds to a different path in its corresponding construction. Undirected paths represent a directed path and its reverse. Paths of length 1 are not included in the visualizations.}
    \label{fig:n=5:case2}
\end{figure}


\end{proof}\bigskip

Theorems \ref{Thm:BalDCTS:Nec} and \ref{Thm:BNHDPD:Suff5} give the complete spectrum for BNHDPD$(5,k)$ since all triples $x_1,x_2,x_3$ that satisfy Theorem \ref{Thm:BalDCTS:Nec} also satisfy Theorem \ref{Thm:BalDCTS:Nec2}. We summarize this result in Theorem \ref{Thm:BalDCTS:Spectrum5}.

\begin{theorem}\label{Thm:BalDCTS:Spectrum5}
    A BNHDPD$(5,k)$ exists if and only if the following conditions all hold:
    \begin{enumerate}
        \item $x_1+2x_2+3x_3=20$
        \item $5\ |\ x_1+x_2+x_3$
    \end{enumerate}
\end{theorem}

\subsection{Sufficiency for BNHDPD$(6,k)$}\label{ssec:BNHDPD:n=6}

In this subsection, we establish sufficient conditions for BNHDPD$(n,k)$ with $n=6$. This is equivalent to constructing BDCTS$(n)$ that contain 6 classes of simple functions and that assess the Chain Rule in isolation.

\begin{theorem}\label{Thm:BNHDPD:Suff6}
    A BNHDPD$(6,k)$ comprised of $x_1$ copies of $\mathcal{D}P_2$, $x_2$ copies of $\mathcal{D}P_3$, $x_3$ copies of $\mathcal{D}P_4$, and $x_4$ copies of $\mathcal{D}P_5$ exists for all nonnegative integer solutions of $x_1+2x_2+3x_3+4x_4=30$ such that $6\ |\ x_1+x_2+x_3+x_4$ and $6\ |\ x_2+2x_3+3x_4$.
\end{theorem}

\begin{proof}
    Let $V(\mathcal{D}K_6)=\{v_1,\hdots, v_6\}$. We proceed by construction via case analysis. Since $6\ |\ x_1+x_2+x_3+x_4$, $6\ |\ x_2+2x_3+3x_4$, and $x_1+x_2+x_3+x_4\leq 30$, we only have to verify solutions of $x_1+x_2+x_3+x_4=s$ for $s\in\{12,18,24\}$. The case $s=30$ is already established by Lemma \ref{Lem:BNHDPD:Trivial}. Within each case, we denote each subcase by the ordered 4-tuple $(x_1,x_2,x_3,x_4)$ of the corresponding solution to $x_1+x_2+x_3+x_4=s$, for $s\in\{12,18,24\}$. For each subcase, we construct the BNHDPD$(6,k)$
    \[\bigcup_{i=1}^{n-2} (P_i\cup P'_i).\] where $P_i, 1\leq i\leq n-2$, is the set of directed paths of length $i$, and  $P'_i, 1\leq i\leq n-2$, is the set of some, but not necessarily all, reverses of $P_i$. Visualizations of each construction are given in Figures \ref{fig:n=6:case1}, \ref{fig:n=6:case2} and \ref{fig:n=6:case3}.
    \medskip

    \noindent\textbf{Case 1:} $x_1+x_2+x_3+x_4=12$
    \begin{enumerate}[label=\textit{Subcase 1\alph*.}]
    
       {\setlength\itemindent{38pt} \item $(3,0,9,0)$}
       
       $P_3=\{v_4v_6v_1v_2, v_6v_2v_3v_4, v_2v_4v_5v_6, v_3v_2v_6v_5, v_1v_6v_4v_3, v_5v_4v_2v_1, v_4v_1v_5v_3, v_5v_1v_3v_6,$\\       
       $ v_6v_3v_5v_2\}$

       $P_1=\{v_1v_4, v_2v_5, v_3v_1\}$

       $P_4,P'_4,P'_3,P_2,P'_2,P'_1=\varnothing$

       {\setlength\itemindent{38pt} \item $(2,2,8,0)$}
       
        $P_3=\{v_4v_6v_1v_2, v_6v_2v_3v_4, v_2v_4v_5v_6, v_2v_5v_3v_6\}$

        $P'_3=\{v_2v_1v_6v_4, v_4v_3v_2v_6, v_6v_5v_4v_2, v_6v_3v_5v_2\}$
        
        $P_2=\{v_5v_1v_3\}$

        $P'_2=\{v_3v_1v_5\}$

        $P_1=\{v_1v_4\}$

        $P'_1=\{v_4v_1\}$

        $P_4,P'_4=\varnothing$

       {\setlength\itemindent{38pt} \item $(1,4,7,0)$}
       
        $P_3=\{v_4v_6v_1v_2, v_6v_2v_3v_4, v_2v_4v_5v_6, v_5v_4v_2v_1, v_1v_6v_4v_3, v_3v_2v_6v_5, v_1v_3v_5v_2\}$
        
        $P_2=\{v_2v_5v_1, v_4v_1v_5, v_5v_3v_6, v_3v_1v_4\}$

        $P_1=\{v_6v_3\}$

        $P_4,P'_4, P'_3,P'_2,P'_1=\varnothing$

       {\setlength\itemindent{38pt} \item $(0,6,6,0)$}
       
       $P_3=\{v_4v_6v_1v_2, v_6v_2v_3v_4, v_2v_4v_5v_6\}$

       $P'_3=\{v_2v_1v_6v_4, v_4v_3v_2v_6, v_6v_5v_4v_2\}$
       
       $P_2=\{v_1v_3v_6, v_3v_5v_2, v_5v_1v_4\}$

       $P'_2=\{v_6v_3v_1, v_2v_5v_3, v_4v_1v_5\}$

       $P_4,P'_4,P_1,P'_1=\varnothing$

       {\setlength\itemindent{38pt} \item $(4,1,4,3)$}
       
        $P_4=\{v_1v_2v_4v_5v_3, v_2v_3v_5v_6v_4, v_3v_4v_6v_1v_5\}$
        
        $P_3=\{v_2v_1v_6v_3, v_4v_3v_2v_6, v_6v_5v_4v_1, v_5v_1v_3v_6\}$
        
        $P_2=\{v_6v_2v_5\}$

        $P_1=\{v_1v_4, v_3v_1, v_4v_2, v_5v_2\}$
       
        $P'_4,P'_3,P'_2,P'_1=\varnothing$

       {\setlength\itemindent{38pt} \item $(3,3,3,3)$}
       
        $P_4=\{v_1v_2v_4v_5v_3, v_2v_3v_5v_6v_4, v_3v_4v_6v_1v_5\}$
        
        $P_3=\{v_2v_1v_6v_3, v_4v_3v_2v_6, v_6v_5v_4v_1\}$
        
        $P_2=\{v_6v_2v_5, v_1v_3v_6, v_5v_1v_4\}$

        $P_1=\{v_3v_1, v_5v_2, v_4v_2\}$

        $P'_4,P'_3,P'_2,P'_1=\varnothing$
       
       \newpage

       {\setlength\itemindent{38pt} \item $(2,5,2,3)$}
       
        $P_4=\{v_1v_2v_4v_5v_3, v_2v_3v_5v_6v_4, v_3v_4v_6v_1v_5\}$
        
        $P_3=\{v_2v_1v_6v_3, v_6v_5v_4v_1\}$
        
        $P_2=\{v_6v_2v_5, v_1v_3v_6, v_5v_1v_4,
        v_4v_3v_2, v_4v_2v_6\}$

        $P_1=\{v_5v_2, v_3v_1\}$
       
        $P'_4,P'_3,P'_2,P'_1=\varnothing$

       {\setlength\itemindent{38pt} \item $(1,7,1,3)$}
       
        $P_4=\{v_1v_2v_4v_5v_3, v_2v_3v_5v_6v_4, v_3v_4v_6v_1v_5\}$
        
        $P_3=\{v_2v_1v_6v_3\}$
        
        $P_2=\{v_6v_2v_5, v_1v_3v_6, v_5v_1v_4,
        v_4v_3v_2, v_4v_2v_6, v_5v_4v_1, v_6v_5v_2\}$

        $P_1=\{v_3v_1\}$

        $P'_4,P'_3,P'_2,P'_1=\varnothing$

       {\setlength\itemindent{38pt} \item $(0,9,0,3)$}
       
        $P_4=\{v_1v_2v_3v_5v_4, v_3v_4v_5v_1v_6, v_5v_6v_1v_3v_2\}$
        
        $P_2=\{v_2v_1v_5, v_6v_5v_3, v_4v_3v_1,
        v_5v_2v_4, v_3v_6v_2,
         v_1v_4v_6\}$

        $P_2=\{v_4v_2v_5,
        v_2v_6v_3, v_6v_4v_1\}$
       
        $P'_4,P_3,P'_3,P_1,P'_1=\varnothing$

       \end{enumerate}\bigskip

\begin{figure}
    \centering
    \input{fig_n=6_Case1}
    \caption{Visual Constructions from Theorem \ref{Thm:BNHDPD:Suff6}, Case 1}
    \floatfoot{Each color corresponds to a different path in its corresponding construction. Undirected paths represent a directed path and its reverse. Paths of length 1 are not included in the visualizations.}
    \label{fig:n=6:case1}
\end{figure}

    \noindent\textbf{Case 2:} $x_1+x_2+x_3+x_4=18$
    \begin{enumerate}[label=\textit{Subcase 2\alph*.}]
    
       {\setlength\itemindent{38pt} \item $(12,0,6,0)$}
       
        $P_3=\{v_1v_2v_4v_5, v_2v_1v_5v_4,
        v_5v_6v_2v_3, v_6v_5v_3v_2,
        v_3v_4v_6v_1, v_4v_3v_1v_6\}$

        $P_1=\{v_2v_6, v_4v_2, v_6v_4, v_1v_3, v_5v_1, v_3v_5, v_1v_4, v_2v_5, v_3v_6\}$

        $P'_1=\{v_4v_1, v_5v_2, v_6v_3\}$
       
        $P_4,P'_4,P'_3,P_2,P'_2=\varnothing$

       {\setlength\itemindent{38pt} \item $(11,2,5,0)$}
       
        $P_3=\{v_1v_2v_4v_5, v_2v_1v_5v_4,
        v_3v_4v_6v_1, v_4v_3v_1v_6,
        v_5v_6v_2v_3\}$
        
        $P_2=\{v_6v_5v_3, v_6v_3v_2\}$

        $P_1=\{v_2v_6, v_4v_2, v_6v_4, v_1v_3, v_5v_1, v_3v_5,v_1v_4, v_2v_5,v_6v_3\}$

        $P'_1=\{v_4v_1, v_5v_2\}$

        $P_4,P'_4,P'_3,P'_2=\varnothing$

       {\setlength\itemindent{38pt} \item $(10,4,4,0)$}
       
        $P_3=\{v_1v_2v_4v_5, v_2v_1v_5v_4, v_3v_4v_6v_1, v_5v_6v_2v_3\}$
        
        $P_2=\{v_6v_5v_3, v_6v_3v_2,
        v_4v_3v_1, v_4v_1v_6\}$

        $P_1=\{v_2v_6, v_4v_2, v_6v_4, v_1v_3, v_5v_1, v_3v_5,v_1v_4, v_2v_5,v_6v_3\}$

        $P'_1=\{v_5v_2\}$

        $P_4,P'_4,P'_3,P'_2=\varnothing$

       {\setlength\itemindent{38pt} \item $(9,6,3,0)$}
       
        $P_3=\{v_3v_4v_6v_1, v_5v_6v_2v_3,
        v_1v_2v_4v_5\}$
        
        $P_2=\{v_6v_5v_3, v_6v_3v_2,
        v_4v_3v_1, v_4v_1v_6,
        v_2v_1v_5, v_2v_5v_4\}$

        $P_1=\{v_2v_6, v_4v_2, v_6v_4, v_1v_3, v_5v_1, v_3v_5,v_1v_4, v_2v_5,v_6v_3\}$

        $P_4,P'_4,P'_3,P'_2,P'_1=\varnothing$

       {\setlength\itemindent{38pt} \item $(8,8,2,0)$}
       
        $P_3=\{v_3v_4v_6v_1, v_5v_6v_2v_3\}$
        
        $P_2=\{v_6v_5v_3, v_6v_3v_2,
        v_4v_3v_1, v_4v_1v_6,
        v_2v_1v_5, v_2v_5v_4,
        v_1v_2v_4, v_1v_4v_5\}$

        $P_1=\{v_1v_3, v_2v_6, v_3v_5, v_3v_6, v_4v_2, v_5v_2, v_5v_1, v_6v_4\}$

        $P_4,P'_4,P'_3,P'_2,P'_1=\varnothing$

       \newpage

       {\setlength\itemindent{38pt} \item $(7,10,1,0)$}
       
        $P_3=\{v_5v_6v_2v_3\}$
        
        $P_2=\{v_6v_5v_3, v_6v_3v_2,
        v_4v_3v_1, v_4v_1v_6,
        v_2v_1v_5, v_2v_5v_4,
        v_1v_2v_4, v_1v_4v_5,
        v_3v_4v_6, v_3v_6v_1\}$

        $P_1=\{v_1v_3, v_2v_6, v_3v_5, v_4v_2, v_5v_2, v_5v_1, v_6v_4\}$

        $P_4,P'_4,P'_3,P'_2,P'_1=\varnothing$

       {\setlength\itemindent{38pt} \item $(6,12,0,0)$}
       
        $P_2=\{v_6v_5v_3, v_6v_3v_2,
        v_4v_3v_1, v_4v_1v_6,
        v_2v_1v_5, v_2v_5v_4,
        v_1v_2v_4, v_1v_4v_5,
        v_3v_4v_6, v_3v_6v_1,$\\
        $v_5v_6v_2, v_5v_2v_3\}$

        $P_1=\{v_1v_3, v_2v_6, v_3v_5, v_4v_2, v_5v_1, v_6v_4\}$

        $P_4,P'_4,P_3,P'_3,P'_2,P'_1=\varnothing$

\end{enumerate}\bigskip

\begin{figure}
    \centering
    \begin{tabular}{|c|c|c|}
\hline
& & \\

\begin{tikzpicture}
\Vertex[x=0,y=2,color=white,label=$v_1$]{1}
\Vertex[x=1.732,y=1,color=white,label=$v_2$]{2}
\Vertex[x=1.732,y=-1,color=white,label=$v_3$]{3}
\Vertex[x=0,y=-2,color=white,label=$v_4$]{4}
\Vertex[x=-1.732,y=-1,color=white,label=$v_5$]{5}
\Vertex[x=-1.732,y=1,color=white,label=$v_6$]{6}

\Edge[color=blue,Direct,bend=7](1)(2)
\Edge[color=blue,Direct,bend=7](2)(4)
\Edge[color=blue,Direct,bend=7](4)(5)

\Edge[color=red,Direct,bend=7](2)(1)
\Edge[color=red,Direct,bend=7](1)(5)
\Edge[color=red,Direct,bend=7](5)(4)

\Edge[color=yellow,Direct,bend=7](5)(6)
\Edge[color=yellow,Direct,bend=7](6)(2)
\Edge[color=yellow,Direct,bend=7](2)(3)

\Edge[color=green,Direct,bend=7](6)(5)
\Edge[color=green,Direct,bend=7](5)(3)
\Edge[color=green,Direct,bend=7](3)(2)

\Edge[color=cyan,Direct,bend=7](3)(4)
\Edge[color=cyan,Direct,bend=7](4)(6)
\Edge[color=cyan,Direct,bend=7](6)(1)

\Edge[color=orange,Direct,bend=7](4)(3)
\Edge[color=orange,Direct,bend=7](3)(1)
\Edge[color=orange,Direct,bend=7](1)(6)

\end{tikzpicture} &

\begin{tikzpicture}
\Vertex[x=0,y=2,color=white,label=$v_1$]{1}
\Vertex[x=1.732,y=1,color=white,label=$v_2$]{2}
\Vertex[x=1.732,y=-1,color=white,label=$v_3$]{3}
\Vertex[x=0,y=-2,color=white,label=$v_4$]{4}
\Vertex[x=-1.732,y=-1,color=white,label=$v_5$]{5}
\Vertex[x=-1.732,y=1,color=white,label=$v_6$]{6}

\Edge[color=blue,Direct,bend=7](1)(2)
\Edge[color=blue,Direct,bend=7](2)(4)
\Edge[color=blue,Direct,bend=7](4)(5)

\Edge[color=red,Direct,bend=7](2)(1)
\Edge[color=red,Direct,bend=7](1)(5)
\Edge[color=red,Direct,bend=7](5)(4)

\Edge[color=yellow,Direct,bend=7](5)(6)
\Edge[color=yellow,Direct,bend=7](6)(2)
\Edge[color=yellow,Direct,bend=7](2)(3)

\Edge[color=cyan,Direct,bend=7](3)(4)
\Edge[color=cyan,Direct,bend=7](4)(6)
\Edge[color=cyan,Direct,bend=7](6)(1)

\Edge[color=orange,Direct,bend=7](4)(3)
\Edge[color=orange,Direct,bend=7](3)(1)
\Edge[color=orange,Direct,bend=7](1)(6)

\Edge[color=pink,Direct,bend=7](6)(5)
\Edge[color=pink,Direct,bend=7](5)(3)

\Edge[color=teal,Direct,bend=7](6)(3)
\Edge[color=teal,Direct,bend=7](3)(2)

\end{tikzpicture} &

\begin{tikzpicture}
\Vertex[x=0,y=2,color=white,label=$v_1$]{1}
\Vertex[x=1.732,y=1,color=white,label=$v_2$]{2}
\Vertex[x=1.732,y=-1,color=white,label=$v_3$]{3}
\Vertex[x=0,y=-2,color=white,label=$v_4$]{4}
\Vertex[x=-1.732,y=-1,color=white,label=$v_5$]{5}
\Vertex[x=-1.732,y=1,color=white,label=$v_6$]{6}

\Edge[color=blue,Direct,bend=7](1)(2)
\Edge[color=blue,Direct,bend=7](2)(4)
\Edge[color=blue,Direct,bend=7](4)(5)

\Edge[color=red,Direct,bend=7](2)(1)
\Edge[color=red,Direct,bend=7](1)(5)
\Edge[color=red,Direct,bend=7](5)(4)

\Edge[color=yellow,Direct,bend=7](5)(6)
\Edge[color=yellow,Direct,bend=7](6)(2)
\Edge[color=yellow,Direct,bend=7](2)(3)

\Edge[color=cyan,Direct,bend=7](3)(4)
\Edge[color=cyan,Direct,bend=7](4)(6)
\Edge[color=cyan,Direct,bend=7](6)(1)

\Edge[color=pink,Direct,bend=7](6)(5)
\Edge[color=pink,Direct,bend=7](5)(3)

\Edge[color=teal,Direct,bend=7](6)(3)
\Edge[color=teal,Direct,bend=7](3)(2)

\Edge[color=violet,Direct,bend=7](4)(3)
\Edge[color=violet,Direct,bend=7](3)(1)

\Edge[color=WildStrawberry,Direct,bend=7](4)(1)
\Edge[color=WildStrawberry,Direct,bend=7](1)(6)

\end{tikzpicture} \\

Subcase (a) & Subcase (b) & Subcase (c) \\
\hline

& &  \\

\begin{tikzpicture}
\Vertex[x=0,y=2,color=white,label=$v_1$]{1}
\Vertex[x=1.732,y=1,color=white,label=$v_2$]{2}
\Vertex[x=1.732,y=-1,color=white,label=$v_3$]{3}
\Vertex[x=0,y=-2,color=white,label=$v_4$]{4}
\Vertex[x=-1.732,y=-1,color=white,label=$v_5$]{5}
\Vertex[x=-1.732,y=1,color=white,label=$v_6$]{6}

\Edge[color=blue,Direct,bend=7](1)(2)
\Edge[color=blue,Direct,bend=7](2)(4)
\Edge[color=blue,Direct,bend=7](4)(5)

\Edge[color=yellow,Direct,bend=7](5)(6)
\Edge[color=yellow,Direct,bend=7](6)(2)
\Edge[color=yellow,Direct,bend=7](2)(3)

\Edge[color=cyan,Direct,bend=7](3)(4)
\Edge[color=cyan,Direct,bend=7](4)(6)
\Edge[color=cyan,Direct,bend=7](6)(1)

\Edge[color=pink,Direct,bend=7](6)(5)
\Edge[color=pink,Direct,bend=7](5)(3)

\Edge[color=teal,Direct,bend=7](6)(3)
\Edge[color=teal,Direct,bend=7](3)(2)

\Edge[color=violet,Direct,bend=7](4)(3)
\Edge[color=violet,Direct,bend=7](3)(1)

\Edge[color=WildStrawberry,Direct,bend=7](4)(1)
\Edge[color=WildStrawberry,Direct,bend=7](1)(6)

\Edge[color=BurntOrange,Direct,bend=7](2)(1)
\Edge[color=BurntOrange,Direct,bend=7](1)(5)

\Edge[color=BlueViolet,Direct,bend=7](2)(5)
\Edge[color=BlueViolet,Direct,bend=7](5)(4)

\end{tikzpicture} &

\begin{tikzpicture}
\Vertex[x=0,y=2,color=white,label=$v_1$]{1}
\Vertex[x=1.732,y=1,color=white,label=$v_2$]{2}
\Vertex[x=1.732,y=-1,color=white,label=$v_3$]{3}
\Vertex[x=0,y=-2,color=white,label=$v_4$]{4}
\Vertex[x=-1.732,y=-1,color=white,label=$v_5$]{5}
\Vertex[x=-1.732,y=1,color=white,label=$v_6$]{6}

\Edge[color=yellow,Direct,bend=7](5)(6)
\Edge[color=yellow,Direct,bend=7](6)(2)
\Edge[color=yellow,Direct,bend=7](2)(3)

\Edge[color=cyan,Direct,bend=7](3)(4)
\Edge[color=cyan,Direct,bend=7](4)(6)
\Edge[color=cyan,Direct,bend=7](6)(1)

\Edge[color=pink,Direct,bend=7](6)(5)
\Edge[color=pink,Direct,bend=7](5)(3)

\Edge[color=teal,Direct,bend=7](6)(3)
\Edge[color=teal,Direct,bend=7](3)(2)

\Edge[color=violet,Direct,bend=7](4)(3)
\Edge[color=violet,Direct,bend=7](3)(1)

\Edge[color=WildStrawberry,Direct,bend=7](4)(1)
\Edge[color=WildStrawberry,Direct,bend=7](1)(6)

\Edge[color=BurntOrange,Direct,bend=7](2)(1)
\Edge[color=BurntOrange,Direct,bend=7](1)(5)

\Edge[color=BlueViolet,Direct,bend=7](2)(5)
\Edge[color=BlueViolet,Direct,bend=7](5)(4)

\Edge[color=CadetBlue,Direct,bend=7](1)(2)
\Edge[color=CadetBlue,Direct,bend=7](2)(4)

\Edge[color=Orchid,Direct,bend=7](1)(4)
\Edge[color=Orchid,Direct,bend=7](4)(5)

\end{tikzpicture} &

\begin{tikzpicture}
\Vertex[x=0,y=2,color=white,label=$v_1$]{1}
\Vertex[x=1.732,y=1,color=white,label=$v_2$]{2}
\Vertex[x=1.732,y=-1,color=white,label=$v_3$]{3}
\Vertex[x=0,y=-2,color=white,label=$v_4$]{4}
\Vertex[x=-1.732,y=-1,color=white,label=$v_5$]{5}
\Vertex[x=-1.732,y=1,color=white,label=$v_6$]{6}

\Edge[color=yellow,Direct,bend=7](5)(6)
\Edge[color=yellow,Direct,bend=7](6)(2)
\Edge[color=yellow,Direct,bend=7](2)(3)

\Edge[color=pink,Direct,bend=7](6)(5)
\Edge[color=pink,Direct,bend=7](5)(3)

\Edge[color=teal,Direct,bend=7](6)(3)
\Edge[color=teal,Direct,bend=7](3)(2)

\Edge[color=violet,Direct,bend=7](4)(3)
\Edge[color=violet,Direct,bend=7](3)(1)

\Edge[color=WildStrawberry,Direct,bend=7](4)(1)
\Edge[color=WildStrawberry,Direct,bend=7](1)(6)

\Edge[color=BurntOrange,Direct,bend=7](2)(1)
\Edge[color=BurntOrange,Direct,bend=7](1)(5)

\Edge[color=BlueViolet,Direct,bend=7](2)(5)
\Edge[color=BlueViolet,Direct,bend=7](5)(4)

\Edge[color=CadetBlue,Direct,bend=7](1)(2)
\Edge[color=CadetBlue,Direct,bend=7](2)(4)

\Edge[color=Orchid,Direct,bend=7](1)(4)
\Edge[color=Orchid,Direct,bend=7](4)(5)

\Edge[color=NavyBlue,Direct,bend=7](3)(4)
\Edge[color=NavyBlue,Direct,bend=7](4)(6)

\Edge[color=Gray,Direct,bend=7](3)(6)
\Edge[color=Gray,Direct,bend=7](6)(1)

\end{tikzpicture}\\

Subcase (d) & Subcase (e) & Subcase (f)\\

\hline

& &  \\

\begin{tikzpicture}
\Vertex[x=0,y=2,color=white,label=$v_1$]{1}
\Vertex[x=1.732,y=1,color=white,label=$v_2$]{2}
\Vertex[x=1.732,y=-1,color=white,label=$v_3$]{3}
\Vertex[x=0,y=-2,color=white,label=$v_4$]{4}
\Vertex[x=-1.732,y=-1,color=white,label=$v_5$]{5}
\Vertex[x=-1.732,y=1,color=white,label=$v_6$]{6}

\Edge[color=pink,Direct,bend=7](6)(5)
\Edge[color=pink,Direct,bend=7](5)(3)

\Edge[color=teal,Direct,bend=7](6)(3)
\Edge[color=teal,Direct,bend=7](3)(2)

\Edge[color=violet,Direct,bend=7](4)(3)
\Edge[color=violet,Direct,bend=7](3)(1)

\Edge[color=WildStrawberry,Direct,bend=7](4)(1)
\Edge[color=WildStrawberry,Direct,bend=7](1)(6)

\Edge[color=BurntOrange,Direct,bend=7](2)(1)
\Edge[color=BurntOrange,Direct,bend=7](1)(5)

\Edge[color=BlueViolet,Direct,bend=7](2)(5)
\Edge[color=BlueViolet,Direct,bend=7](5)(4)

\Edge[color=CadetBlue,Direct,bend=7](1)(2)
\Edge[color=CadetBlue,Direct,bend=7](2)(4)

\Edge[color=Orchid,Direct,bend=7](1)(4)
\Edge[color=Orchid,Direct,bend=7](4)(5)

\Edge[color=NavyBlue,Direct,bend=7](3)(4)
\Edge[color=NavyBlue,Direct,bend=7](4)(6)

\Edge[color=Gray,Direct,bend=7](3)(6)
\Edge[color=Gray,Direct,bend=7](6)(1)

\Edge[color=blue,Direct,bend=7](5)(6)
\Edge[color=blue,Direct,bend=7](6)(2)

\Edge[color=red,Direct,bend=7](5)(2)
\Edge[color=red,Direct,bend=7](2)(3)

\end{tikzpicture} & & \\

Subcase (g) &  & \\

\hline

\end{tabular}
    \caption{Visual Constructions from Theorem \ref{Thm:BNHDPD:Suff6}, Case 2}
    \floatfoot{Each color corresponds to a different path in its corresponding construction. Undirected paths represent a directed path and its reverse. Paths of length 1 are not included in the visualizations.}
    \label{fig:n=6:case2}
\end{figure}

    \noindent\textbf{Case 3:} $x_1+x_2+x_3+x_4=24$
    \begin{enumerate}[label=\textit{Subcase 3\alph*.}]
    
       {\setlength\itemindent{38pt} \item $(21,0,3,0)$}
       
      $P_3=\{v_6v_1v_2v_4, v_2v_3v_4v_6, v_4v_5v_6v_2\}$

      $P_1\cup P'_1=E(\mathcal{D}K_6)\setminus P_3$

      $P_4,P'_4,P'_3,P_2,P'_2=\varnothing$

       {\setlength\itemindent{38pt} \item $(20,2,2,0)$}
       
      $P_3=\{v_2v_3v_4v_6, v_4v_5v_6v_2\}$
      
      $P_2=\{v_6v_1v_3, v_1v_2v_4\}$
      
      $P_1\cup P'_1=E(\mathcal{D}K_6)\setminus (P_3\cup P_2)$

      $P_4,P'_4,P'_3,P'_2=\varnothing$

       {\setlength\itemindent{38pt} \item $(19,4,1,0)$}
       
      $P_3=\{v_4v_5v_6v_2\}$
      
      $P_2=\{v_1v_2v_4, v_2v_3v_5, v_3v_4v_6, v_6v_1v_3\}$

      $P_1\cup P'_1=E(\mathcal{D}K_6)\setminus (P_3\cup P_2)$

      $P_4,P'_4,P'_3,P'_2=\varnothing$

       {\setlength\itemindent{38pt} \item $(18,6,0,0)$}
       
      $P_2=\{v_1v_2v_4, v_2v_3v_5,
      v_3v_4v_6,
      v_4v_5v_1,
      v_5v_6v_2,
      v_6v_1v_3\}$

      $P_1\cup P'_1=E(\mathcal{D}K_6)\setminus P_2$

      $P_4,P'_4,P_3,P'_3,P'_2=\varnothing$
      
        \end{enumerate} \end{proof}

Theorems \ref{Thm:BalDCTS:Nec}, \ref{Thm:BalDCTS:Nec2}, and \ref{Thm:BNHDPD:Suff6} give the complete spectrum for BNHDPD$(6,k)$. We summarize this result in Theorem \ref{Thm:BalDCTS:Spectrum6}.

\begin{theorem}\label{Thm:BalDCTS:Spectrum6}
    A BNHDPD$(6,k)$ exists if and only if the following conditions all hold:
    \begin{enumerate}
        \item $x_1+2x_2+3x_3+4x_4=30$
        \item $6\ |\ x_1+x_2+x_3+x_4$
        \item $6\ |\ x_2+2x_3+3x_4$
    \end{enumerate}  
\end{theorem}

\begin{figure}
    \centering
    \begin{tabular}{|c|c|c|c|}
\hline
& & & \\

\begin{tikzpicture}[scale=0.75]
\Vertex[x=0,y=2,color=white,label=$v_1$]{1}
\Vertex[x=1.732,y=1,color=white,label=$v_2$]{2}
\Vertex[x=1.732,y=-1,color=white,label=$v_3$]{3}
\Vertex[x=0,y=-2,color=white,label=$v_4$]{4}
\Vertex[x=-1.732,y=-1,color=white,label=$v_5$]{5}
\Vertex[x=-1.732,y=1,color=white,label=$v_6$]{6}

\Edge[color=blue,Direct](6)(1)
\Edge[color=blue,Direct](1)(2)
\Edge[color=blue,Direct](2)(4)

\Edge[color=red,Direct](2)(3)
\Edge[color=red,Direct](3)(4)
\Edge[color=red,Direct](4)(6)

\Edge[color=yellow,Direct](4)(5)
\Edge[color=yellow,Direct](5)(6)
\Edge[color=yellow,Direct](6)(2)

\end{tikzpicture} &

\begin{tikzpicture}[scale=0.75]
\Vertex[x=0,y=2,color=white,label=$v_1$]{1}
\Vertex[x=1.732,y=1,color=white,label=$v_2$]{2}
\Vertex[x=1.732,y=-1,color=white,label=$v_3$]{3}
\Vertex[x=0,y=-2,color=white,label=$v_4$]{4}
\Vertex[x=-1.732,y=-1,color=white,label=$v_5$]{5}
\Vertex[x=-1.732,y=1,color=white,label=$v_6$]{6}

\Edge[color=blue,Direct](6)(1)
\Edge[color=blue,Direct](1)(3)

\Edge[color=green,Direct](1)(2)
\Edge[color=green,Direct](2)(4)

\Edge[color=red,Direct](2)(3)
\Edge[color=red,Direct](3)(4)
\Edge[color=red,Direct](4)(6)

\Edge[color=yellow,Direct](4)(5)
\Edge[color=yellow,Direct](5)(6)
\Edge[color=yellow,Direct](6)(2)

\end{tikzpicture} &

\begin{tikzpicture}[scale=0.75]
\Vertex[x=0,y=2,color=white,label=$v_1$]{1}
\Vertex[x=1.732,y=1,color=white,label=$v_2$]{2}
\Vertex[x=1.732,y=-1,color=white,label=$v_3$]{3}
\Vertex[x=0,y=-2,color=white,label=$v_4$]{4}
\Vertex[x=-1.732,y=-1,color=white,label=$v_5$]{5}
\Vertex[x=-1.732,y=1,color=white,label=$v_6$]{6}

\Edge[color=blue,Direct](6)(1)
\Edge[color=blue,Direct](1)(3)

\Edge[color=green,Direct](1)(2)
\Edge[color=green,Direct](2)(4)

\Edge[color=red,Direct](2)(3)
\Edge[color=red,Direct](3)(5)

\Edge[color=cyan,Direct](3)(4)
\Edge[color=cyan,Direct](4)(6)

\Edge[color=yellow,Direct](4)(5)
\Edge[color=yellow,Direct](5)(6)
\Edge[color=yellow,Direct](6)(2)

\end{tikzpicture} &

\begin{tikzpicture}[scale=0.75]
\Vertex[x=0,y=2,color=white,label=$v_1$]{1}
\Vertex[x=1.732,y=1,color=white,label=$v_2$]{2}
\Vertex[x=1.732,y=-1,color=white,label=$v_3$]{3}
\Vertex[x=0,y=-2,color=white,label=$v_4$]{4}
\Vertex[x=-1.732,y=-1,color=white,label=$v_5$]{5}
\Vertex[x=-1.732,y=1,color=white,label=$v_6$]{6}

\Edge[color=blue,Direct](6)(1)
\Edge[color=blue,Direct](1)(3)

\Edge[color=green,Direct](1)(2)
\Edge[color=green,Direct](2)(4)

\Edge[color=red,Direct](2)(3)
\Edge[color=red,Direct](3)(5)

\Edge[color=cyan,Direct](3)(4)
\Edge[color=cyan,Direct](4)(6)

\Edge[color=yellow,Direct](4)(5)
\Edge[color=yellow,Direct](5)(1)

\Edge[color=orange,Direct](5)(6)
\Edge[color=orange,Direct](6)(2)

\end{tikzpicture}\\

Subcase (a) & Subcase (b) & Subcase (c) & Subcase (d) \\

\hline

\end{tabular}
    \caption{Visual Constructions from Theorem \ref{Thm:BNHDPD:Suff6}, Case 3}
    \floatfoot{Each color corresponds to a different path in its corresponding construction. Undirected paths represent a directed path and its reverse. Paths of length 1 are not included in the visualizations.}
    \label{fig:n=6:case3}
\end{figure}
    

\section{Sample \emph{Designed} Task Sets}\label{sec:HWsamples}

The task sets given in this section are constructed using Subcase 1d of the proof of Theorem \ref{Thm:BNHDPD:Suff5}. First, we show a task set where the skills are derivatives of specific functions. \newpage

\noindent\textbf{Task Set 1.} Let $1\mapsto x^2$, $2\mapsto \sin x$, $3\mapsto \ln x$, $4\mapsto e^x$, $5\mapsto\arctan x$. Then, the balanced task set is
\begin{multicols}{2}
\begin{enumerate}
    \item $e^{\sin x}$
    \item $\sin({\arctan x})$
    \item $\ln(\arctan x)$
    \item $\arctan(\ln(x^2))$
    \item $\arctan(\sin(e^x))$
    \item $\displaystyle e^{(\ln x)^2}$
    \item $(\sin(\ln x))^2$
    \item $\displaystyle(e^{\ln(\sin x)})^2$
    \item $\displaystyle\ln(e^{\arctan(x^2)})$
    \item $\sin((\arctan(e^x))^2)$
\end{enumerate}
\end{multicols}

Next, we show a task set where the skills are derivatives of simple function classes and values within each class are randomly chosen.  \medskip

\noindent\textbf{Task Set 2.} Let $1\in\{x^n\ |\ 2\leq n\leq10\}$, $2\in\{\sin x,\cos x,\tan x,\sec x,\csc x,\cot x\}$, $3\in\{\ln x,\log_a x\ \text{where}\ 2\leq a\leq10\}$, $4\in\{e^x, a^x\ \text{where}\ 2\leq a\leq 10\}$, $5\in\{\arcsin x,\arccos x,\arctan x\}$. To randomly choose values for each task, the elements of each set defined above were converted to integer values, which were then chosen using the following Excel 2020 formulas.
\begin{enumerate}
    \item \textbf{Power functions} were chosen using \texttt{=RANDBETWEEN(2,10)}, where the value of this formula corresponds to the exponent of the power function.
    \item \textbf{Trigonometric functions} were chosen using \texttt{=RANDBETWEEN(1,6)}, where the value of this function corresponds to one of the six trigonometric functions as ordered above.
    \item \textbf{Exponential and logarithmic functions} were chosen using a two-step formula. First, \texttt{=RANDBETWEEN(0,1)} chose between a natural exponential or logarithm and a general exponential or logarithm. This was followed by \texttt{=IF(}[cell containing 0 or 1]\texttt{=1,RANDBETWEEN(2,10),0)} to choose the base of a general exponential or logarithm.
    \item \textbf{Inverse trigonometric functions} were chosen using \texttt{=RANDBETWEEN(1,3)}, where the value of this function corresponds to one of the three inverse trigonometric functions as ordered above. 
\end{enumerate}\newpage

The first balanced task set produced by this code was
\begin{multicols}{2}
\begin{enumerate}
    \item $7^{\sec x}$
    \item $\csc({\arccos x})$
    \item $\log_8(\arccos x)$
    \item $\arctan(\ln(x^7))$
    \item $\arcsin(\csc(e^x))$
    \item $\displaystyle 6^{(\log_7 x)^8}$
    \item $(\csc(\log_5 x))^5$
    \item $\displaystyle(e^{\ln(\csc x)})^4$
    \item $\displaystyle\ln(3^{\arcsin(x^{10})})$
    \item $\csc((\arccos(e^x))^{10})$
\end{enumerate}
\end{multicols}

\section{Conclusion}\label{sec:conc}

In this paper, we established existence theorems for balanced non-Hamiltonian directed path decompositions of the complete directed graph. This extends prior work by Meszka \& Skupie\'n \cite{Meszka} by requiring that the decomposition be balanced. Furthermore, these results are applicable to the construction of task sets in Calculus I that assess the Chain Rule due to a correspondence between labeled directed paths and composite functions. 

\subsection{Future Directions}

A natural future avenue to explore is the existence of balanced (or near-balanced) task sets assessing additional derivative computation skills: sums, products, and quotients of simple functions of one variable. Since these are binary operations, we will need to decompose the complete directed graph into subgraphs containing directed stars. Preliminary mathematical results by Colbourn, Hoffman, and Rodger \cite{DSD} are established for the existence of directed star decompositions of the complete directed graph. 

Alternatively, we could consider a different type of practice in the task set. Rather than constructing mixed practice task sets, as we did here with the complete directed graph, we can consider the complete directed multigraph with index $\lambda$ for task sets with a hybrid mixed-block practice. Meszka \& Skupie\'n's work on non-Hamiltonian directed paths \cite{Meszka} and Colbourn, Hoffman, \& Rodger's work on directed stars \cite{DSD} contain analogous results regarding decompositions of the complete directed multigraph with index $\lambda$.

\subsection{Closing Remarks}

The applications of this work demonstrate that while Combinatorial Designs are a highly abstract area of discrete mathematics, they have the potential for meaningful applications in undergraduate mathematics education. Since such applications are, as of now, widely understudied, it invites the possibility of rich interdisciplinary collaboration, uniting two fields, Combinatorial Designs and Undergraduate Mathematics Education, in novel ways.

\end{document}